\documentclass[a4paper,12pt]{article} 
\usepackage[left=2.5cm,right=2.5cm,top=2.5cm,bottom=2.5cm,bindingoffset=0cm]{geometry}
\usepackage{cmap}
\usepackage[T2A,T2B]{fontenc}
\usepackage[cp1251]{inputenc} 
\usepackage[english]{babel} 
\usepackage{indentfirst}
\usepackage{amssymb,amsmath}
\usepackage{pst-all}
\usepackage{graphicx}
\usepackage[hyphens]{url} 
\usepackage{algpseudocode} 
\usepackage[off]{auto-pst-pdf}
\usepackage{amsthm}
\ifpdf
  \usepackage{hyperref}
  \usepackage[hyphens]{url} 
\fi

\usepackage{color,ifpdf,graphicx}
\ifpdf 
  \DeclareGraphicsExtensions{.pdf,.png,.mps}
\else 
  \DeclareGraphicsExtensions{.eps,.bmp}
  \usepackage{pstricks}
\fi
\usepackage{pgf,epic,bez123}

\title{Numerical homotopy continuation for control and online identification of nonlinear systems: the survey of selected results}

\date{}
\author{Alex Borisevich\footnotemark[1]}

\newcommand{\Lie}{{\mathcal{L}}}
\newcommand{\R}{\mathbb{R}}
\DeclareMathOperator{\const}{const}
\DeclareMathOperator{\rank}{rank}

\newtheorem{Def}{Definition}
\newtheorem{Thm}{Theorem}

\newtheorem{Exmpl}{Example}
\newenvironment{Prf}[1][Proof]{\begin{trivlist}
\item[\hskip \labelsep {\bfseries #1}]}{\end{trivlist}}
\theoremstyle{remark}
\newtheorem{Rem}{Remark}

\begin{document}

\maketitle

\ifpdf
\footnotetext[1]{St. Petersburg State Polytechnical University, Polytekhnicheskaya 29, St. Petersburg, 195251, Russia, \href{mailto:alex.borysevych@spbstu.ru}{alex.borysevych@spbstu.ru}}
\fi

\begin{abstract}
The article gives an overview of the parameter numerical continuation methodology applied to setpoint control and parameter identification of nonlinear systems. The control problems for affine systems as well as general (nonaffine) nonlinear systems are considered. Online parameter identification is also presented in two versions: with linear and nonlinear nonconvex parameterization. Simulation results for illustrative examples are shown.
\end{abstract}

\section*{Introduction}

The theory of nonlinear systems is an active area of research in which exist well-established paradigms and in the same time there are many important and difficult questions unanswered. If the problem of analysis (stability, sensitivity to parameter variations) can be resolved at a sufficient level for practical use, including by means of simulation (MATLAB, Ansys), the problem of controller synthesis for nonlinear systems are still "hot spot" in control theory.

It is difficult to list all of the approaches to the control of nonlinear systems, but if we try to name at least the most famous of them, there are:

\begin{itemize}

\item Local linearization near working point and application of the linear robust control to nonlinear systems \cite{1}. Considering nonlinearity as uncertainty of the nonlinear model, the robust regulator synthesized on the basis of linearization can be used for control of nonlinear system.

\item Feedback linearization \cite{2,3,4}. In this approach nonlinear dynamics of the initial system is transformed to the linear by means of nonlinear feedback. Then obtained linear system can be controlled with tools of the linear control theory.

\item Control on the basis of differential flatness \cite{5}. If the states and inputs of systems can be parametrized via the output variables and their derivatives, such a system is called the differentially flat. This allows to explicitly solve the problem of the inversion of the system and generate the control actions.

\item The use of control Lyapunov functions or their generalizations \cite{6}. Synthesis of Lyapunov functions for the specified control can be formulated as an optimization problem of semi-definite programming \cite{7}. Special method for a particular class of systems management is a integrator backstepping \cite{8}. This method uses a recursive decomposition of the object to the subsystems of a special type, which are stabilized by means of Lyapunov functions.

\item Nonlinear model predictive control \cite{9}. 

\item Variable structure control, which also includes a relay and a sliding mode control \cite{10}. Promising direction is the approximation of nonlinear system by the set of hybrid systems and application of hybrid control \cite{11}.

\item Intelligent control based on neural networks or fuzzy logic.

\end{itemize}

In this paper, we will lean on the feedback linearization approach. This method is attractive because it allows to implement controls for nonlinear systems as simple as for linear. Constructing of linearization transformation is based on the analytical differentiation of the system outputs and can easily be done manually or automatically \cite{12}.

From the practical point of view, a class of objects that are suitable for control by feedback linearization is substantially less than, for example, for the method of Lyapunov function or sliding mode control. The main difficulty is the mandatory requirement for constant order of state model equations, known as the constancy of the output relative degree. When relative degree is ill-defined then linearizing coordinate transformation is singular in some regions of phase space and the control action is discontinuous.

Problems with the relative degree of dynamic systems arise naturally in underactuated mechanical systems, where the number of degrees of freedom is greater than installed actuators \cite{13}. Region of phase space in which the relative degree of the system changed is related to the configurations in which the system changes its controllability.

This paper provides an overview of the new approach to control of nonlinear systems, which combines the conceptual simplicity of feedback linearization methods and at the same time extends the applicability to irregular systems with ill-defined output relative degree.

Approach described below is based on numerical parameter continuation methods for solving systems of nonlinear equations \cite{19}. This method is based on parameterized combination of the original problem and some simple problem with a known solution. The immediate motivation for the use of parameter continuation methods in control is a series of papers \cite{14,15}, which addresses the application of these methods directly in the physical experiments.

The article consists of several parts. The first part defines the basic concepts like models of nonlinear systems, control problems and methods. The second part is devoted the following objectives: output setpoint control and identification of parameters in real time (online). After a description of the method a simple example is given that covers a particular method.

The material is based on previous papers \cite{16,17,18}, but all of the methods described below are new and not previously published.

\section{Preliminaries}

\subsection{Plant models}

In this paper we consider nonlinear MIMO systems with $m$ inputs and $m$ outputs in the state space of dimension $n$:

\begin{equation}\label{plant}
\dot{x} = f(x, u), \; y = h(x)
\end{equation}

where $x \in X \subseteq \R^n$, $y \in Y \subseteq \R^m$, $u \in U \subseteq \R^m$, and maps $f : \R^n \times \R^m \to \R^n$, $h : \R^n \to \R^m$ are smooth vector fields $f,h \in C^{\infty}$. The functions $f(.)$ and $h(.)$ are considered bounded in $X$.

One important particular case of \eqref{plant} is affine nonlinear system of form 

\begin{equation}\label{plant_affine}
\dot{x} = f(x) + \sum_{i=1}^m g_i(x) u_i, \; y = h(x)
\end{equation}

where the functions $g_i : \R^n \to \R^n$ are bounded smooth on $X$ vector field.

Also for the problems of adaptive control and parameter identification considered the following control objects with parametric uncertainty

\begin{equation}\label{plant_uncertain}
\dot{x} = f(x, \theta, u), \; y = h(x, \theta)
\end{equation}

\begin{equation}\label{plant_affine_uncertain}
\dot{x} = f(x, \theta) + \sum_{i=1}^m g_i(x, \theta) u_i, \; y = h(x, \theta)
\end{equation}

where $\theta \in \Theta \subseteq \R^n$ is an unknown constant parameter which is not available for measuring or monitoring.

\subsection{Control problems}

In this paper we consider the following control problems for nonlinear systems of form \eqref{plant} or \eqref{plant_affine}:

\begin{itemize}

\item Output setpoint control for the system \eqref{plant_affine} in affine form.

\item Output setpoint control for the general nonlinear system \eqref{plant}.

\item Online parameter identification of nonlinear affine \eqref{plant_affine_uncertain} and nonaffine \eqref{plant_uncertain} plants.

\end{itemize}

Let us give formal definitions of the problems considered in this paper \cite{23}:

\begin{Def}
The control problem of output setpoint tracking is to synthesize of such control law $u(t) = u(x)$, which asymptotically drives the system output $y$ of system \eqref{plant} to constant value $\bar y = \const$: $\lim_{t\to\infty} y(t) = \bar y$. In the particular case $\bar y = 0$ the control is called output zeroing problem.
\end{Def}

Without loss of generality, everywhere below we assume that $\bar y = 0$ and the objective of control is output zeroing.

\begin{Def}
The on-line identification problem for unknown constant parameter $\hat \theta$ is to find an adaptation algorithm for estimate $\hat \theta(t)$, which would minimize the difference $y^*(t) - y(t)$ between the measured system response $y^*(t)$ to the input $u(t)$ and the response of model \eqref{plant_uncertain} to the same input. $\square$
\end{Def}

\subsection{Linear approximation in point}\label{subsec:lin_point}

Let there be given object of control in the form of \eqref{plant} with initial conditions $x(0) = x_0$, $u(0) = u_0$. The task is to stabilize \eqref{plant} locally near the initial point $x_0 \in X$, $u_0 \in U$. Effective method is the linear approximation of \eqref{plant} in a small neighborhood around $(x_0, u_0)$ and dropping the high-order nonlinear terms of the power expansion \cite{23}:

\begin{equation}\label{lin_appox}
\begin{gathered}
\dot x = f(x,u) \approx f(x_0,u_0) + D_x(x_0,u_0)(x - x_0) + D_u(x_0,u_0)(u - u_0) = A x + B u + \Delta_x \\
y = h(x) \approx h(x_0) + D_h(x_0)(x - x_0) = C x + \Delta_y
\end{gathered}
\end{equation}

where $D_x f(x,u) = \left ( \dfrac{\partial f_i}{\partial x_j}\right )$, $D_u f(x,u) = \left ( \dfrac{\partial f_i}{\partial u_j} \right )$, $D_x h(x) = \left ( \dfrac{\partial h_i}{\partial x_j} \right )$.

\subsection{Feedback linearization}\label{subsec:lin_feedback}

\begin{Def} 
The MIMO affine nonlinear system \eqref{plant_affine} has relative degree $r_j$ for output $y_j$ in $\mathfrak{S} \subseteq \R^n$ if at least for one function $g_i$ is true

\begin{equation}
\Lie_{g_i} \Lie_f^{r_j-1} h_j \ne 0
\end{equation}

where $\Lie_f \lambda = \dfrac{\partial \lambda(x)}{\partial x} f(x) = \displaystyle\sum\nolimits_{i=1}^n \dfrac{\partial \lambda(x)}{\partial x_i} f_i(x)$ is a Lie derivative of function $\lambda$ along a vector field $f$. 
\end{Def} 

It means that at least one input $u_k$ influences to output $y_j$ after $r_j$ differentiations.

Number $r = \sum_{i=1}^m r_j$ is called as the total relative degree of system. If $ r = n $ and matrix 

\begin{equation}\label{decmatrix}
A(x) = \begin{pmatrix} \Lie_{g_1} \Lie_f^{r_1-1} h_1(x) & \cdots & \Lie_{g_m} \Lie_f^{r_1-1} h_1(x) \\ \vdots & \cdots & \vdots \\ \Lie_{g_1} \Lie_f^{r_m-1} h_m(x) & \cdots & \Lie_{g_m} \Lie_f^{r_m-1} h_m(x) \end{pmatrix}
\end{equation}

is full rank, then the original dynamical system \eqref{plant_affine} in $\mathfrak{S}$ equivalent to system:

\begin{equation}\label{affine_out}
y^{(r_j)}_j = \Lie_f^{r_j} h_j + \sum_{i=1}^m \Lie_{g_i} \Lie_f^{r_j-1} h_j \cdot u_i = B(x) + A(x) \cdot u
\end{equation}

The nonlinear feedback

\begin{equation}\label{nltrans}
u = A(x)^{-1} [ v - B(x) ]
\end{equation}

converts in subspace $\mathfrak{S}$ original dynamical system \eqref{plant_affine} to linear:

\begin{equation}\label{lin}
y^{(r_j)} = v_j
\end{equation}

\subsection{Linearization of irregular systems}

\begin{Def}
The relative degree of affine system \eqref{plant_affine} is called ill-defined, if there are points $x^{\circ} \in \mathfrak{S}$ in the state space, in which the decoupling matrix $A(x)$ is rank deficient, i.e.
\begin{equation}
\det A(x^{\circ}) = 0
\end{equation}
\end{Def}

Systems with ill-defined relative degree is also called irregular.

Irregularity does not mean that it is impossible to solve the particular control problem (output setpoint control or parameter identification). It just means that not all the trajectories of the linear system \eqref{lin} to which transforms the original nonlinear system \eqref{plant_affine} can be reproduced at the output $y$ due to the nature of nonlinearities.

\begin{Exmpl}
The simplest specific example is the system $\dot x = u$, $y = h(x) = x(x^2 - 1) + 1$, $x(0) = 1$ for which the problem of output zeroing $y \to 0$ is needed to solve.

If the system under consideration was a constant relative degree, the use of control $v = -y$ after feedback linearization would give the output trajectory of $ y (t) = \exp (-t) $, which is everywhere decreasing $\dot y(t) < 0$.

In this case, the nonlinearity $y = h(x)$ has two limit points $x^{\circ}_{1,2} = \pm 1/\sqrt{3}$, in which $h'_x(x^{\circ}_{1,2}) = 0$. Any trajectory $y(t)$, that connects $y(0) = 1$ with $y(T) = 0$ passes sequentially through the points $y^{\circ}_{1} = h(3^{-1/2})$ and $y^{\circ}_{2} = h(- 3^{-1/2})$, and besides $y^{\circ}_{2} > y^{\circ}_{1}$. Hence, any trajectory $y(t)$ on the interval $(0, t_1)$ should decrease with time (Figure 1), on the interval $(t_1, t_2)$ increase, and in the interval again decrease. Such a trajectory is not reproducible using the feedback linearization.
\end{Exmpl}

\subsection{Numerical continuation method}

Let it is necessary to solve system of the nonlinear equations

\begin{equation}\label{nleq}
\phi(\xi) = 0
\end{equation}

where $\phi: \R^m \to \R^m$ is vector-valued smooth nonlinear function.

Lets $\Omega \subset \R^m$ is open set and $C(\bar \Omega)$ is set of continuous maps from its closure $\bar \Omega$ to $\R^m$. Functions $F_0, F_1 \in C(\bar \Omega)$ are homotopic (homotopy equivalent) if there exists a continuous mapping

\begin{equation}\label{homotopy}
H : \bar \Omega \times [0,1] \to \R^m
\end{equation}

that $H(\xi,0) = F_0(\xi)$ and $H(\xi,1) = F_1(\xi)$ for all $\xi \in \bar \Omega$. It can be shown \cite{9} that the equation $H(\xi, \lambda) = 0$ has solution $(\xi, \lambda)$ for all $\lambda \in [0,1]$. The objective of all numerical continuation methods is tracing of implicitly defined function $H(\xi, \lambda) = 0$ for $\lambda \in [0,1]$.

Lets $H: \mathcal{D} \to \R^m$ is $C^1$-continuous function on an open set $\mathcal{D} \subset \R^{m+1}$, and the Jacobian matrix $DH(\xi,\lambda)$ is full-rank $\operatorname{rank}DH(\xi,\lambda) = m$ for all $(\xi,\lambda) \in \mathcal{D}$. Then, for all $(\xi,\lambda) \in \mathcal{D}$ exists a unique vector $\tau \in \R^{m + 1}$ such as

\begin{equation}
DH(\xi,\lambda) \cdot \tau = 0, \;
\|\tau\|_2 = 1, \;
\det \begin{pmatrix} DH(\xi,\lambda) \\ \tau^T \end{pmatrix} > 0,
\end{equation}

and mapping

\begin{equation}\label{tanfunc}
\Psi: \mathcal{D} \to \R^{m+1}, \; \Psi: (\xi,\lambda) \mapsto \tau
\end{equation}

is locally Lipschitz on $\mathcal{D}$.

Function \eqref{tanfunc} specifies the autonomous differential equation

\begin{equation}\label{homode}
\frac{d}{dt} \begin{pmatrix} \xi \\ \lambda \end{pmatrix} = \Psi(\xi,\lambda), \; \xi(0) = \xi_0, \; \lambda(0) = 0
\end{equation}

which has a unique solution $\xi(t),\lambda(t)$ according to a theorem of solution existence for the Cauchy problem.

It can be shown \cite{9} that integral curve $\gamma(t) = (\xi(t),\lambda(t))$ reaches the solution point $(\xi^{*}, 1)$ where $\phi(\xi^{*}) = 0$ in finite time $t < \infty$.

\subsection{Control of the output as a two-step problem}\label{subsec:two_stage}.

The following fact is important for the methods discussed in subsections \ref{sec:affine} and \ref{sec:nonaffine}.

Let there is a system of the form \eqref{plant}. Required to solve the output setpoint control problem $y \to \bar y$. Let $\bar x = h^{-1} (\bar y)$. Then the problem can be divided into two phases:

\begin{itemize}

\item Implementation of movement along the trajectory $x(t)$, starting at $x(0) = x_0$ and ending in $x(T) = x_e$, where $x_e \in \epsilon(\bar x)$ is point in the neighborhood of $\bar x$, $T < \infty$.

\item asymptotic stabilization \eqref{plant} near the point $\bar x$ for $t > T$ with initial conditions $x(T) = x_e$.

\end{itemize}

The main result of the methods in sections \ref{sec:affine} and \ref{sec:nonaffine} is explicitly regard to task 1. The last task 2 may be solved by the linearization at $\bar x$ according to section \ref{subsec:lin_point}.

\section{Main result}

\subsection{Regularization of affine nonlinear system by parameter continuation}\label{sec:affine}

Let's associate with the plant \eqref{plant_affine} linear dynamics system with $m$ inputs $u$, $n$ states $z$, $m$ outputs $\eta$ and with the same relative degries $r_i$ for outputs such as in \eqref{plant_affine}

\begin{equation}\label{linsys}
\dot z = A z + B u, \; \eta = C z, \; \frac{d^{(r_i)}}{d t^{(r_i)}}\eta = u_i.
\end{equation}

Now we want to mix output of \eqref{plant_affine} with output of \eqref{linsys}

\begin{equation}\label{H}
H = (1 - \lambda) \cdot \eta + \lambda \cdot y
\end{equation}

where $\lambda \in [0,1]$ is a continuous time-dependant parameter, value of which determines the relative contribution of the outputs $\eta$ and $y$ to $H$.

Control process can be represented as follows. At time $t = 0$ parameter value is $\lambda = 0$, which corresponds to control of only a linear system \eqref{linsys}. After, the value of $\lambda$ starts increasing, but in general, it is described by some bounded function $\lambda(t)$. The parameter $\lambda$ reaches $\lambda = 1$ in a finite time $t = T < \lambda$, which corresponds to the end of the process control.

By definition of the output relative degree, each component $H_i$ in \eqref{H} should be differentiated $r_i$ times with respect to $t$ until it becomes an explicit function of any input $u$. We obtain after differentiation:

\begin{equation}\label{difform_H_i_pre}
\begin{gathered}
H_i^{(r_i)} = - \sum_{k=1}^{r_i-1} C_{r_i}^k \eta_i^{(r_i-k)} \lambda^{(k)} + (1 - \lambda) u_i + (y_i - \eta_i) \lambda^{(r_i)} + \sum_{k=1}^{r_i-1} C_{r_i}^k y_i^{(r_i-k)} \lambda^{(k)} + \\
\quad + \lambda \left ( \Lie_f^{r_i} h_i + \sum_{k=1}^m \Lie_{g_k} \Lie_f^{r_i-1} h_i \cdot u_k \right ) = 0,
\end{gathered}
\end{equation}

where $C_n^k$ are the binomial coefficients.

Let's denote following

\begin{equation}\label{A_B_H}
\begin{gathered}
\Lambda_i = (\lambda, \dot \lambda, \ddot \lambda, ... , \lambda^{(r_i-1)}) \\
\mathcal{A}_{i,1}(x,z,\Lambda_i) = \left ( \lambda \Lie_{g_k} \Lie_f^{r_i-1} h_i \right )^T + (1 - \lambda) \cdot \delta_k \\
\mathcal{A}_{i,2}(x,z) = y_i - \eta_i \\
\mathcal{B}_i(x,z,\Lambda_i) = - \sum_{k=1}^{r_i-1} C_{r_i}^k \eta_i^{(r_i-k)} \lambda^{(k)} + \sum_{k=1}^{r_i-1} C_{r_i}^k y_i^{(r_i-k)} \lambda^{(k)} + \lambda \Lie_f^{r_i} h_i 
\end{gathered}
\end{equation}

where $\delta_k$ is a vector of dimension $m$, all of whose components are equal to 0 except the $k$-th.

With introduction of notation from \eqref{A_B_H} the equation \eqref{difform_H_i_pre} can be written more compactly:

\begin{equation}\label{difform_H_i}
H^{(r_i)} = \mathcal{A}_{i,1}(x,z,\Lambda) \cdot u + \mathcal{A}_{i,2}(x,z) \cdot \lambda^{(r_i)} + \mathcal{B}_i(x,z,\Lambda),
\end{equation}

Considering all of the components $H_i$ after differentiation according to the relative degrees of outputs $r_i$ it is possible to write an algebraic condition \eqref{difform_H_i} in the matrix form:

\begin{equation}\label{H_all}
\begin{split}
H^{(r)} &= \mathcal{A}_{1}(x,z,\Lambda) \cdot u + \mathcal{A}_{2}(x,z) \cdot \lambda^{(r_{max})} + \mathcal{B}(x,z,\Lambda) \\
 &= \mathcal{A}(x,z,\Lambda) \cdot \begin{pmatrix} u \\ \lambda^{(r_{max})} \end{pmatrix} + \mathcal{B}(x,z,\Lambda) \\
\end{split}
\end{equation}

where $r_{max} = \max\{ r_i \}$, $\Lambda = (\lambda, \dot \lambda, \ddot \lambda, ... , \lambda^{(r_{max}-1)})$, $H^{(r)} = (H_1^{(r_1)}, H_2^{(r_2)}, ... H_{m}^{(r_m)})^T$.

Now we ready to formulate following theorem, which is our main constructive result

\begin{Thm}\label{Thm:sol}

Suppose

\begin{equation}\label{ExistCond}
\operatorname{rank}\mathcal{A}(x,z,\Lambda) = m,
\end{equation}

then for state feedback with bounded new input $v$

\begin{equation}\label{Sol}
\begin{pmatrix} u \\ \lambda^{(r_{max})} \end{pmatrix} = \alpha \cdot \tau(\mathcal{A}) + \mathcal{A}^+ (v - \mathcal{B}), \\
\end{equation}

where vector $\tau(\mathcal{A})$ calculated to satisfy following constraints 

\begin{equation}\label{tau_calc}
\mathcal{A} \cdot \tau = 0, \; \|\tau\|_2 = 1, \; \det \begin{pmatrix} \mathcal{A} \\ \tau^T \end{pmatrix} > 0
\end{equation}

and $\alpha = \operatorname{const} \in \R_{+}$, $z(0) = 0$, $\Lambda(0) = 0$ following can be stated:

0. There exists smooth control trajectory $(u(t),\lambda(t))$ generated by \eqref{Sol} which leaves point $(u(0),0)$.

1. System \eqref{H_all} transformed by feedback \eqref{Sol} to linear controllable form 

\begin{equation}\label{H_lin}
H^{(r)} = v
\end{equation}

2. Curve $(u(t), \lambda(t))$ either passes through point $\lambda = 1$ or diffeomorphic to a circle.

\end{Thm}

\begin{Prf}
Full proof is omitted. Assertion 0 is proved by explicit expression of vector $\tau(\mathcal{A})$ for $t = 0$. Proposition 1 is trivial by substituting \eqref{Sol} into \eqref{H_all}. Statement 2 is argued as analogous statement for the solution curve in parameter continuation method for nonlinear equations \cite{19}.
\end{Prf}

\begin{Rem}

As a result of feedback \eqref{Sol} system with the output \eqref{H} is transformed into a linear \eqref{H_lin}, which can be controlled by any known methods of linear control theory (pole placement, linear quadratic regulator, etc.). Outer linear feedback is crucial for the practical applications of the proposed approach, as it allows to compensate uncertainties in the model and implement robust control systems.

\end{Rem}

\begin{Rem}

The last theorem indicates that the parameter $\alpha$ is additional degree of freedom in the controller design. The larger this constant, the faster the solution arrives to the $\lambda = 1$, but numerical integration becomes more stiff.

With a decrease of the parameter $\alpha$ to less than a certain critical value $\alpha_{min} $ convergence of parameter $\lambda$ to the value of $\lambda(t) \to \lambda^* < 1$ is observed, which is unacceptable. Accurate characterization of this phenomenon, as well as assessment of $\alpha_{min}$, is the subject of further research.

\end{Rem}

\begin{Rem}

In some cases, it is desirable to vary the speed of convergence $\lambda(t)$ to the final value of the parameter $\lambda(t) \to \lambda^* = 1$ for the regulation of transient time. Let us introduce the scalar coefficient $\gamma > 0$, with which possible to control speed of change $\lambda(t)$. If $\gamma > 1$ then $\lambda(t)$ reaches a value of $\lambda^* = 1$ faster than compared to \eqref{Sol}. Otherwise, when $\gamma < 1$ transition process is slowing down. 

We modify the \eqref{Sol} as follows 

\begin{equation}\label{Sol_gamma}
\begin{pmatrix} u \\ \lambda^{(r_{max})} \end{pmatrix} = \alpha \cdot \gamma \cdot \tau(\mathcal{A}) + Q (\mathcal{A} \cdot Q)^+ (v - \mathcal{B}),
\end{equation}

where $Q = diag(1,...1,\gamma)$ is square diagonal matrix, the last diagonal element of which is equal to $\gamma$.

Equation \eqref{Sol_gamma} obtained from \eqref{Sol} by scaling vector $\tau$ and application of weighted by matrix $Q$ pseudoinverse of $A$.

\end{Rem}

\begin{Rem}

Attention is required in case $r_i > 1$ for the relative degree of the output $H$. Fulfillment of condition $H^{(r)} = 0$ when $r_i > 1$, in general, does not consequent to $H = 0$. For example, if $\ddot H(t) = 0$ for a time interval $t \in [0, T]$, but if initially $\dot H(0) = h_0 \ne 0$, then for the whole time interval is obtained $\dot H(t) = h_0$ and $H(t) = h_0 \cdot t \ne 0$.

In the case of SISO system with the choice of initial conditions for the $\lambda^{(i)}(0) $ is possible to achieve $H^{(i)}(0) = 0 $. For MIMO systems with $r_i > 1$ the selection of $\lambda^{(i)}(0)$ in general can only minimize $H_k^{(i)}$. Hence, necessary to use of an external linear feedback to control the $H$ to 0 in case of $r_i > 1$.

\end{Rem}

\begin{Rem}

Condition \eqref{ExistCond} is a standard assumption of using parameter continuation method, which corresponds to the possible existence of limit points of trajectories $(u(t), \lambda(t))$ at which $\mathcal{A}_1 \notin \operatorname{im} \mathcal{A}_2$, and the absence of bifurcation points. At the same time in some regions of phase space $X \times Z$ may be a situation where $\operatorname{rank} \mathcal{A}_1(x,z,\Lambda) < m$, in that case, the system with output \eqref{H} cannot be directly linearized by the feedback, but the proposed method is still applicable.

\end{Rem}

\begin{Rem}

Condition \eqref{ExistCond} can be relaxed a little, but we do not consider this here. In fact, the proposed method allows the existence of the phase space of simple bifurcation points where $\dim \ker \mathcal{A} = 2$. When the control trajectory passes through a simple bifurcation point the sign of vector $\tau$ flipped (a more detailed analysis in \cite{9}).

Overcoming the bifurcation points, in which is observed $\mathcal{A}_1 \in \operatorname{im} \mathcal{A}_2$, also possible within the known approaches for the numerical parameter continuation (e.g., using the Lyapunov-Schmidt decomposition \cite{9}).

\end{Rem}

\subsection{Illustrative example of affine irregular system control}\label{subsect:ex1}

Consider following abstract example of MIMO system, that changes its relative degree in the state space

\begin{equation}\label{plant_ex}
\begin{gathered}
\dot x_1 = u_1 + x_2^3 \\
\dot x_2 = u_2 + x_1^3 \\
y_1 = x_1^3 - x_1 + 1 \\
y_2 = x_2^4 \cos(2 x_2)
\end{gathered}
\end{equation}

with initial conditions $x(0) = (1,1)^T$. We need to solve the problem of output zeroing $y \to 0$. 

Differentiating the outputs, we obtain

\begin{equation}\label{plant_ex_dy}
\begin{gathered}
\dot y_1 = \left (3 x_1^2 - 1 \right ) \cdot (u_1 + x_2^3) = a_{11} u_1 + b_1 \\
\dot y_2 = \left ( 4 x_2^3 \cos(2 x_2) - 2 x_2^4 \sin(2 x_2) \right ) \cdot (u_2 + x_1^3) = a_{22} u_2 + b_2
\end{gathered}
\end{equation}

Obviously, the system in interval $x \in [0,1]^2$ can not be completely linearized by the feedback, because there are exists such $x^*$ that $a_{11}(x^*) = 0$ or $a_{22}(x^*) = 0$.

Let's associate with \eqref{plant_ex} linear system of a form

\begin{equation}\label{plant_ex_lin}
\begin{gathered}
\dot \eta_1 = u_1, \; \dot \eta_2 = u_2
\end{gathered}
\end{equation}

with initial conditions $\eta(0) = (0,0)^T$

According to the equation \eqref{difform_H_i_pre} we obtain for $\dot H = 0$ the following

\begin{equation}\label{H_all_ex_1}
\mathcal{A}_{1} = \lambda \begin{pmatrix} a_{11} & 0 \\ 0 & a_{22} \end{pmatrix} + (1-\lambda) E, \; \mathcal{A}_{2} = y - \eta, \; \mathcal{B} = \lambda \begin{pmatrix} b_1 \\ b_2 \end{pmatrix}
\end{equation}

To test the robustness of the controller to the system output $y$ was applied additive perturbation of the form $\Delta y = (1, \sin(20 t))^T$. External control circuit is implemented using a P-controller with a gain of 100. The value of $\alpha$ is chosen as $ \alpha = 20 $.

The model in Simulink to control the system \label{plant_ex} shown in figure 1. Modeling results are shown on figures 2-4.

\begin{center}
\includegraphics[width=1.0\textwidth]{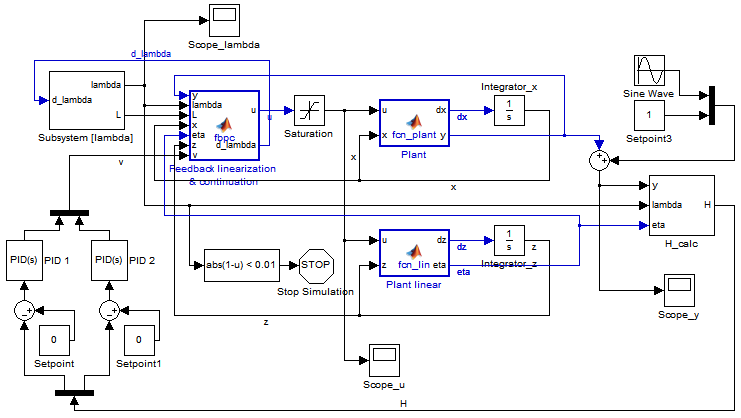}

Figure 1. Simulink model.
\end{center}

\begin{center}
\includegraphics[width=0.85\textwidth]{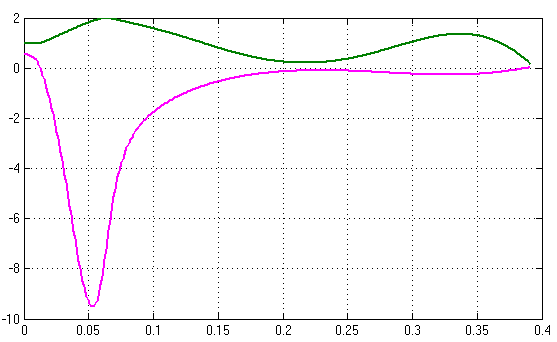}

Figure 2. Output response.
\end{center}

\begin{center}
\includegraphics[width=0.85\textwidth]{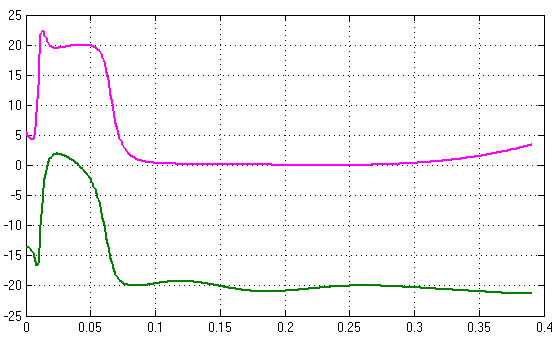}

Figure 3. Input controls.
\end{center}

\begin{center}
\includegraphics[width=0.85\textwidth]{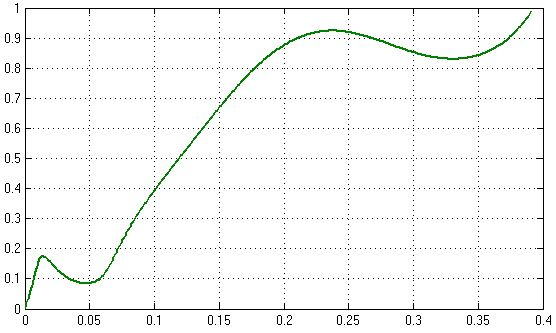}

Figure 4. Dynamics of parameter $\lambda$.
\end{center}

\subsection{Control of nonaffine systems}\label{sec:nonaffine}

Formally, the nonaffine systems of the form \eqref{plant} can be transformed to affine form by adding integrator to inputs: 

\begin{equation}\label{plnat_affinize}
\begin{gathered}
\dot x = f(x, x_u) \\
\dot x_u = u \\
y = h(x) 
\end{gathered}
\end{equation}

As the result there is a new state vector $(x, x_u)$ and the relative degree of each output is increased by one. Consider another way for the case when the increase of the relative degree is undesirable.

Let there be given nonlinear system of the general form \eqref{plant}, for which it is necessary to solve a problem of output zeroing $y \to 0$. 
Denoting $y_0 = h(x_0)$, it is possible to write the following modification of \eqref{H}:

\begin{equation}\label{H_linearize}
H = y + y_0 \lambda - y_0 = 0
\end{equation}

Feature of \eqref{H_linearize} in comparison with \eqref{H} is absence of need for system \eqref{linsys} and linearity of $H$ in the variables $y$ and $\lambda$. Note that the material of section \ref{sec:affine} can also be implemented for the function $Н$ of form \eqref{H_linearize} instead of \eqref{H}.

The next step of approach is the linear approximation of function $f(x, u)$ for variable $u$ around point $u_i$

\begin{equation}\label{plant_linearized}
\dot{x} = f(x, u) \approx f(x, u_i) + D_x f(x ,u_i) (u - u_i) = \hat f^i(x) + \hat g^i(x) u
\end{equation}

where

\begin{equation}
\begin{gathered}
\hat f^i(x) = f(x, u_i) - D_x f(x ,u_i) \cdot u_i \\
\hat g^i(x) = D_x f(x ,u_i)
\end{gathered}
\end{equation}

Thus, the system \eqref{plant} approximated by the affine system in the neighborhood of $u_i$. System behavior of \eqref{plant} in a whole space $U$ can be approximated by the bundle of systems \eqref{plant_linearized}, each of which is a linear approximation of for a particular $u_i \in U$.

Since \eqref{plant_linearized} is the affine system, the concept of relative degree output can be applied.

Differentiating each component $H_j$ according to its relative degree $r_j$, can be obtained using \eqref{affine_out} the following

\begin{equation}\label{dot_H_linearize}
H_j^{(r_j)} = \Lie_{\hat f^i}^{r_j} h_j + \sum_{k=1}^m \Lie_{\hat g_k^i} \Lie_{\hat f^i}^{r_j-1} h_j \cdot u_k + y_j(0) \lambda^{(r_j)}
\end{equation}

Denoting

\begin{equation}
\begin{gathered}
\mathsf{A}^i_{1,j}(x) = \left ( \Lie_{\hat g_k^i} \Lie_{\hat f^i}^{r_j-1} h_j \right )^T \\
\mathsf{A}^i_{2,j} = y_j(0) \\ 
\mathsf{B}^i_j(x) = \Lie_{\hat f^i}^{r_j} h_j
\end{gathered}
\end{equation}

equation \eqref{dot_H_linearize} can be rewritten in compact form

\begin{equation}
H_j^{(r_j)} = \mathsf{B}^i_j(x) + \mathsf{A}^i_{1,j}(x) \cdot u + \mathsf{A}^i_{2,j} \cdot \lambda^{(r_j)}
\end{equation}

Considering all components of $H_j$, it is possible to write \eqref{dot_H_linearize} in matrix form by analogy with \eqref{H_all}:

\begin{equation}\label{H_all_linearize}
\begin{split}
H^{(r)} &= \mathsf{A}_{1}^i(x) \cdot u + \mathsf{A}_{2}^i \cdot \lambda^{(r_{max})} + \mathsf{B}^i(x,\Lambda) \\
 &= \mathsf{A}^i(x) \cdot \begin{pmatrix} u \\ \lambda^{(r_{max})} \end{pmatrix} + \mathsf{B}^i(x,\Lambda) \\
\end{split}
\end{equation}

where $r_{max} = \max\{ r_i \}$, $\Lambda = (\lambda, \dot \lambda, \ddot \lambda, ... , \lambda^{(r_{max}-1)})$.

The following theorem provides a theoretical basis for the possibility of parameter continuation application to control the bundle of systems in form of \eqref{plant_linearized}

\begin{Thm}\label{Thm:infinitesimal}

Suppose that for some moment of time $t_i$ defined $x_i$, $u_i$ and $\lambda_i$ and for moment of time $t_{i+1} = t_i + \Delta t$ quantities $x_{i+1}$, $u_{i+1}$ and $\lambda_{i+1}$ are given, satisfying the equations \eqref{H_linearize} and \eqref{plant}, i.e.

\begin{equation}\label{H_linearize_thm_next}
\begin{gathered}
H = h(x_i) + \lambda_i \cdot y_0 - y_0 = 0 \\
H = h(x_{i+1}) + \lambda_{i+1} \cdot y_0 - y_0 = 0
\end{gathered}
\end{equation}

with

\begin{equation}
\begin{gathered}
\lambda_{i+1} = \lambda_{i} + \Delta \lambda_i \\
x_{i+1} = x_i + \Delta x_i
\end{gathered}
\end{equation}

then $\Delta x_i$ and $\Delta \lambda_i$ are infinitesimals of the same order.

\end{Thm} 

\begin{Prf}
Full proof is omitted due to simplicity. The theorem is proved as follows: nonlinear function $h$ in last equation \eqref{H_linearize_thm_next} expanded in power series around $x_i$ and $\lambda_i$. From the expansion $\Delta\ lambda_i$ and $\Delta x_i$ are expressed. Next limit of the ratio $\Delta \lambda_i / \| \Delta x_i \| $ for $\Delta \lambda_i \to 0$ is composed and its boundedness is shown.
\end{Prf}

The theorem \ref{Thm:infinitesimal} provides an important theoretical conclusion: a small change in parameter $\lambda$ leads to a small change in the $x$ satisfying \eqref{H_linearize}. Hence, after a small change of $\lambda$ the control action $u$, which must be applied to the input of the plant \eqref{plant}, can be obtained from the linearization of the system around the previous value of $\lambda$.

Control method can be described as follows:

0. $\lambda = 0$, $u_i := u_0$, where $u_0$ is solution of $f(x_0,u_0) = 0$.

1. Obtain linearization \eqref{plant_linearized} at point $u_i$ as functions $\hat f^i$ and $\hat g^i$.

2. During time interval $\Delta t$ apply the feedback in form \eqref{Sol} for system with outputs \eqref{H_all_linearize}

\begin{equation}\label{Sol_linearize}
\begin{pmatrix} u \\ \lambda^{(r_{max})} \end{pmatrix} = \alpha \cdot \tau(\mathsf{A}) + (\mathsf{A}^i)^+ \cdot [v - \mathsf{B}^i], \\
\end{equation}

where vector $\tau(\mathsf{A})$ calculated in such a way to satisfy the equations \eqref{tau_calc} for matrix $\mathsf{A}^i$, and external $v$ control input for linear system $H^{(r)} = v$ stabilization.

3. Update the variable $u_i := u$, go to point 1 of algorithm. If $\lambda \ge 1$ then exit.

\subsection{An illustrative example of nonaffine nonlinear system control}\label{subsect:ex2}

Considered an abstract system with one input and one output form

\begin{equation}\label{plant_ex_nonaf}
\begin{gathered}
\dot x = u^3 (x^2 + 1) + e^{-u} \\
y = x (x^2 - 1) + 1
\end{gathered}
\end{equation}

with initial state $x(0) = 1$. Necessary to stabilize the system with output zeroing $y \to 0$. 

Let's write down the output function \eqref{H_linearize} of system \eqref{plant_ex_nonaf}

\begin{equation}\label{plant_ex_nonaf_H_linearize}
\begin{split}
H &= y + y_0 \lambda - y_0 \\
& = x (x^2 - 1) + \lambda
\end{split}
\end{equation}

Next transform the system \eqref{plant_ex_nonaf} to affine form \eqref{plant_linearized} by linear approximation of the first equation at point $u_i$:

\begin{equation}\label{plant_ex_nonaf_linearized}
\begin{gathered}
\dot x = \hat f^i(x) + \hat g^i(x) \\
y = x (x^2 - 1) + 1 \\
\hat f^i(x) = u_i^3 (x^2 + 1) + e^{-u_i} + (e^{-u_i} - 3 u_i^2 (x^2 + 1)) u_i \\
\hat g^i(x) = 3 u_i^2 (x^2 + 1) - e^{-u_i}
\end{gathered}
\end{equation}

Differentiating \eqref{plant_ex_nonaf_H_linearize}, it is possible to write the equations \eqref{dot_H_linearize} and \eqref{H_all_linearize} for the considered example

\begin{equation}\label{plant_ex_nonaf_dot_H}
\begin{split}
\dot H &= \begin{pmatrix} \dfrac{\partial h}{\partial x} \cdot \hat g^i(x) & y_0 \end{pmatrix} \cdot \begin{pmatrix} u \\ \dot \lambda \end{pmatrix} + \frac{\partial h}{\partial x} \cdot \hat f^i(x) \\
       &= \mathsf{A}^i(x) \cdot \begin{pmatrix} u \\ \dot \lambda \end{pmatrix} + \mathsf{B}^i(x)
\end{split}
\end{equation}

where $\dfrac{\partial h}{\partial x} = 3 x^2 - 1$ 

Matrices $\mathsf{A}^i$ and $\mathsf{B}^i$ from \eqref{plant_ex_nonaf_dot_H} are used in the feedback \eqref{Sol_linearize} of control system. Outer control loop is implemented as a proportional controller with a coefficient of 10. The value of $\alpha$ is chosen $ \alpha = 2 $.

Model of the control in the MATLAB/Simulink environment \eqref{plant_ex_nonaf} is shown at Figure 5. The simulation results are presented in Figures 6-8.

\begin{center}
\ifpdf 
  \includegraphics[width=1.0\textwidth]{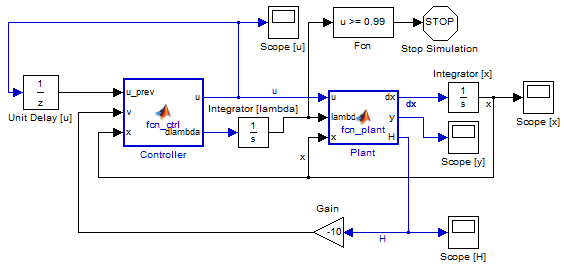}
\fi

Figure 5. Model of closed-loop system in Simulink.
\end{center}

\begin{center}
\ifpdf 
  \includegraphics[width=0.85\textwidth]{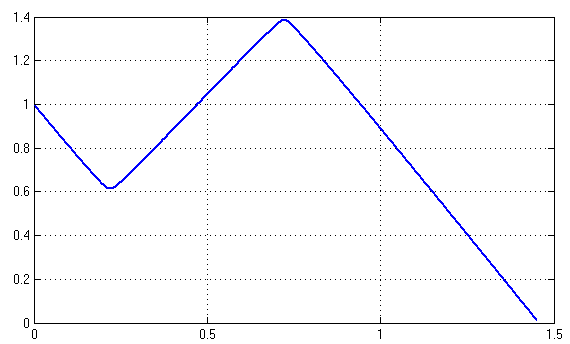}
\fi

Figure 6. Output response $y(t)$.
\end{center}

\begin{center}
\ifpdf 
  \includegraphics[width=0.85\textwidth]{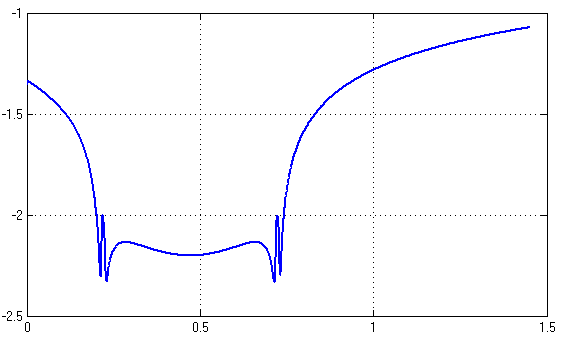}
\fi

Figure 7. Input signals $u(t)$.
\end{center}

\begin{center}
\ifpdf 
  \includegraphics[width=0.85\textwidth]{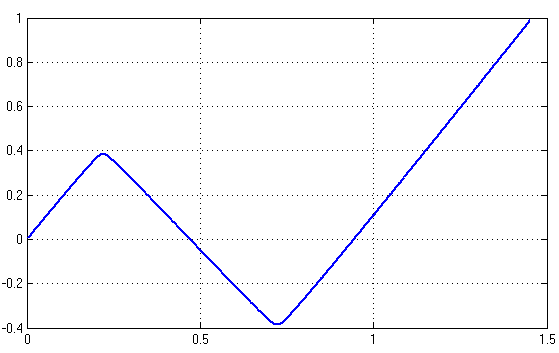}
\fi

Figure 8. Dynamics of parameter $\lambda$.
\end{center}

\subsection{Identification of parameters for system with linear parameterized uncertainty}\label{sec:affine_id}

In this section, the ideas presented in the part \ref{sec:affine}, extended to the problem of parameter identification of nonlinear systems in real time. The purpose of identification is to minimize the difference $x - \hat x$, where $x$ is an output of the process with unknown parameters, and $\hat x$ the output of process model. To solve this problem we introduce a new function $H$ and the way to avoid the digital differentiation
of $x$ is considered.

Let the following nonlinear system in which uncertainties are parametrized linearly by a vector $\theta$ is given

\begin{equation}\label{plant_affine_uncertain_}
\dot{x} = f(x,u) + w(x,u) \cdot \theta,
\end{equation}

where vector $\theta \in \Theta \subset \R^q$ is quasi-constant unknown parameter which should be identified in the process of adaptation, $w : \R^n \times \R^m \to \R^q \times \R^n$ is the know vector function. State $x$ considered available for measurement or observation.

For system \eqref{plant_affine_uncertain_} model is known, which uses a parameter estimate $\hat \theta$ instead of true unknown value $\theta$:

\begin{equation}\label{plant_affine_uncertain_model}
\dot{\hat x} = f(\hat x, u) + w(\hat x, u) \cdot \hat \theta
\end{equation}

with initial conditions $\hat x(0) = \hat x_0$. Next, it is assumed that $\hat x_0 \ne x_0$.

The problem is to find such a value of parameter $\hat \theta$, which would minimize the difference $x - \hat x \to 0$. 

By analogy with \eqref{H_linearize} associate for the systems \eqref{plant_affine_uncertain_} and \eqref{plant_affine_uncertain_model} following function, which takes into account that the final goal of regulation by $\hat \theta$ is zeroing the difference $x - \hat x$. 

\begin{equation}\label{H_linearize_affine_uncertain}
H = x - \hat x - (1-\lambda) (x_0 - \hat x_0)
\end{equation}

Differentiating \eqref{H_linearize_affine_uncertain}, obtain

\begin{equation}\label{dot_H_linearize_affine_uncertain}
\dot H = \dot x - \dot {\hat x} + \dot \lambda (x_0 - \hat x_0)
\end{equation}

equation \eqref{dot_H_linearize_affine_uncertain} can be written using the notation of \eqref{H_all_linearize} for matrices $\mathcal{A}$ and $\mathcal{B}$

\begin{equation}\label{dot_H_linearize_affine_uncertain_AB}
\begin{split}
\dot H &= \dot x - f(\hat x, u) - w(\hat x, u) \cdot \hat \theta + \dot \lambda (x_0 - \hat x_0) \\
&= \mathcal{A}_1 \cdot \hat \theta + \mathcal{A}_2 \cdot \dot \lambda + \mathcal{B}
\end{split}
\end{equation}

where 

\begin{equation}
\begin{gathered}
\mathcal{A}_1 = - w(\hat x, u) \\
\mathcal{A}_2 = x_0 - \hat x_0 \\
\mathcal{B} = \dot x - f(\hat x, u)
\end{gathered}
\end{equation}

Then to \eqref{dot_H_linearize_affine_uncertain_AB} it is possible to apply approach of solution continuation for the equation $\dot H = 0$ by analogy with \eqref{Sol} and \eqref{Sol_linearize}. Hence we obtain the control $\hat \theta(t)$ and dynamics of parameters $\lambda(t)$

\begin{equation}\label{Sol_uncertain}
\begin{pmatrix} \hat \theta \\ \dot \lambda \end{pmatrix} = \alpha \cdot \tau(\mathcal{A}) + \mathcal{A}^+ \cdot (-k H - \mathcal{B})
\end{equation}

where vector $\tau(\mathcal{A})$ calculated in such way as to satisfy the equations \eqref{tau_calc} for $\mathcal{A}$, and $k > 0$ is the proportional loop gain.

A separate task is to calculate $\dot x$, since the variable $x$ it is result of measurement at which there will be noise and which numerical differentiation is undesirable. 

Let's rewrite the last equation in \eqref{Sol_uncertain}

\begin{equation}
\begin{split}
\begin{pmatrix} \hat \theta \\ \dot \lambda \end{pmatrix} &= \alpha \cdot \tau(\mathcal{A}) + \mathcal{A}^+ \cdot (-k H - \mathcal{B}) \\
&= \alpha \cdot \tau(\mathcal{A}) - k \mathcal{A}^+ H - \mathcal{A}^+ \mathcal{B} \\
&= \alpha \cdot \tau(\mathcal{A}) - k \mathcal{A}^+ H + \mathcal{A}^+ f(\hat x, u) - \mathcal{A}^+ \dot x
\end{split}
\end{equation}

Next, we consider the matrix $\mathcal{M}(t) = \alpha \cdot \tau(\mathcal{A}) - k \mathcal{A}^+ H + \mathcal{A}^+ f(\hat x, u)$ and $\mathcal{N}(t) = - \mathcal{A}^+$, then let's split them so that it was possible to write

\begin{equation}\label{Sol_uncertain_for_int}
\begin{split}
\begin{pmatrix} \hat \theta \\ \dot \lambda \end{pmatrix} &= \alpha \cdot \tau - k \mathcal{A}^+ H + \mathcal{A}^+ f(\hat x, u) - \mathcal{A}^+ \dot x \\
&= \mathcal{M} + \mathcal{N} \dot x \\
&= \begin{pmatrix} \mathcal{M}_1 \\ \mathcal{M}_2 \end{pmatrix} + \begin{pmatrix} \mathcal{N}_1 \\ \mathcal{N}_2 \end{pmatrix} \dot x
\end{split}
\end{equation}

Let's integrate $\dot \lambda = \mathcal{M}_2 + \mathcal{N}_2 \dot x$ to obtain the dynamics $\lambda(t)$. Applying integration by parts to the last term of $\mathcal{N}_2 \dot x$, it is possible to write the following equation

\begin{equation}\label{lambda_uncertain_int}
\lambda(t) = \int_0^t \mathcal{M}_2(\tau) d\tau + \mathcal{N}_2(t) \cdot x(t) - \mathcal{N}_2(0) \cdot x(0) - \int_0^t \dot{\mathcal{N}_2}(\tau) \cdot x(\tau) d\tau
\end{equation}

The equation \eqref{lambda_uncertain_int} determines the dynamics of parameter $\lambda(t)$ without direct differentiation of $\dot x$. Unfortunately, a similar equation from \eqref{Sol_uncertain_for_int} for $\theta(t)$ can not be obtained analytically.

Because $\mathcal{N}_2$ calculated on the basis of function $w(\hat x, u)$ from the process model \eqref{plant_affine_uncertain_model} and it does not contain noise, the digital differentiation of $\mathcal{N}_2$ for calculating $ \dot{\mathcal{N}_2}$ applicable in this case.

To calculate the $\hat \theta$ consider a model \eqref{plant_affine_uncertain_model} with outputs \eqref{H_linearize_affine_uncertain} for the purpose of regulation $H \to 0$. Let's perform feedback linearization of system \eqref{plant_affine_uncertain_model} as follows

\begin{equation}\label{plant_affine_uncertain_model_lin}
\hat \theta = [w(\hat x, u)]^+ (v_{\theta} - f(\hat x, u))
\end{equation}

that transforms of nonlinear system \eqref{plant_affine_uncertain_model} to linear $\dot{\hat x} = v_{\theta}$. Further, using proportional control $v_\theta = -k_\theta \cdot H$, $k_\theta > 0$, it is possible to solve the problem of regulating $H$ to $0$. Hence we finally obtain from \eqref{plant_affine_uncertain_model_lin} dynamics of $\hat \theta$

\begin{equation}\label{plant_affine_uncertain_model_lin_ctrl}
\hat \theta = [w(\hat x, u)]^+ (-k_\theta \cdot H - f(\hat x, u))
\end{equation}

Note that the proportional control $v_\theta = -k_\theta \cdot H$ can be complemented with an integral component to improve the accuracy.

Also, it should be noted that, since in the \eqref{plant_affine_uncertain_model_lin_ctrl} the matrix pseudo-inverse of $w(\hat x, u)$ is calculated, which can be degenerate for some pairs $(\hat x, u)$, the trajectory $\hat \theta(t)$, in general, is discontinuous. But as in the output equation \eqref{H_linearize_affine_uncertain} regularizing term introduced $(1-\lambda) (x_0 - \hat x_0)$, it can be expected that the trajectory $\hat \theta$ asymptotically approaches the $\theta$.

\subsection{Illustrative example of identification for linearly parametrized system}\label{subsect:ex3}

The abstract system with one input and one output is considered

\begin{equation}\label{plant_ex_id}
\begin{gathered}
\dot x = (1 - x^2) \theta + u
\end{gathered}
\end{equation}

with the value $\theta = 5$  of the parameter to be identified and initial state $x(0) = x_0 = 0$.

The control action is linearly increasing signal $u(t) = t$.

For the model of the process we will use the following system

\begin{equation}\label{plant_ex_id_model}
\begin{gathered}
\dot{\hat x} = (1 - \hat x^2) \hat \theta + u
\end{gathered}
\end{equation}

with initial state $\hat x(0) = \hat x_0 = 0.5$.

Let's write the function \eqref{H_linearize_affine_uncertain} for a system output \eqref{plant_ex_id}

\begin{equation}\label{plant_ex_id_H}
H = x - \hat x - (1-\lambda) (x_0 - \hat x_0)
\end{equation}

Having differentiated \eqref{plant_ex_id_H}, it is possible to write the equations \eqref{dot_H_linearize_affine_uncertain_AB} for this system

\begin{equation}
\begin{split}
\dot H &= \dot x + (\hat x^2 - 1) \cdot \hat \theta - u + \dot \lambda (x_0 - \hat x_0) \\
&= \mathcal{A}_1 \cdot \hat \theta + \mathcal{A}_2 \cdot \dot \lambda + \mathcal{B}
\end{split}
\end{equation}

where $\mathcal{A}_1 = \hat x^2 - 1$, $\mathcal{A}_2 = x_0 - \hat x_0$, and $\mathcal{B} = \dot x - u$.

Using the notation from \eqref{Sol_uncertain_for_int} equation \eqref{lambda_uncertain_int} acquires the following form:

\begin{equation}\label{plant_ex_id_lambda}
\begin{gathered}
\mathcal{M}(t) = \alpha \cdot \tau  + \begin{pmatrix} \hat x^2 - 1 & x_0 - \hat x_0 \end{pmatrix}^+ (-k H + u) = \begin{pmatrix} \mathcal{M}_1 \\ \mathcal{M}_2 \end{pmatrix} \\
\mathcal{N}(t) = - \begin{pmatrix} \hat x^2 - 1 & x_0 - \hat x_0 \end{pmatrix}^+ = \begin{pmatrix} \mathcal{N}_1 \\ \mathcal{N}_2 \end{pmatrix} \\
\lambda(t) = \int_0^t \mathcal{M}_2(\tau) d\tau + \mathcal{N}_2(t) \cdot x(t) - \mathcal{N}_2(0) \cdot x(0) - \int_0^t \dot{\mathcal{N}_2}(\tau) \cdot x(\tau) d\tau
\end{gathered}
\end{equation}

The dynamics of parameter $\hat \theta$, in this case according to \eqref{plant_affine_uncertain_model_lin_ctrl}

\begin{equation}\label{plant_ex_id_ctrl}
\hat \theta = (1 - \hat x^2)^+ \cdot (-k_\theta \cdot H - u)
\end{equation}

where $(a)^+ = a^{-1}$ for $a \ne 0$, and $(a)^+ = 0$ for $a = 0$.

Model of the control for \eqref{plant_ex_id} in the Simulink environment is shown at Figure 9. The model consists of several subsystems: 

- "Plant" is the plant \eqref{plant_ex_id} itself, 

- "Plant [model] \& Theta" is the model of plant \eqref{plant_ex_id_model} and calculation of parameter $\hat \theta$ dynamics by \eqref{plant_ex_id_ctrl},

- "Lambda calculation" is the block for calculation of $\lambda$ by \eqref{plant_ex_id_lambda}, 

- "Parameter continuation" block for calculation of matrices $\mathcal{M}$ and $\mathcal{N}$.

Feedback constants are chosen as follows: $\alpha = 2.5$, $k = 5$, $k_\theta = 10$. To model the stability of adaptation to measurement noise of state $x$ source of uniformly distributed noise in the range $[-0.01, 0.01]$ is added. Parameter $\hat \theta$ limited in range of $[0, 10]$.

The simulation results are presented in Figures 9-12.

\begin{center}
\ifpdf 
  \includegraphics[width=0.8\textwidth]{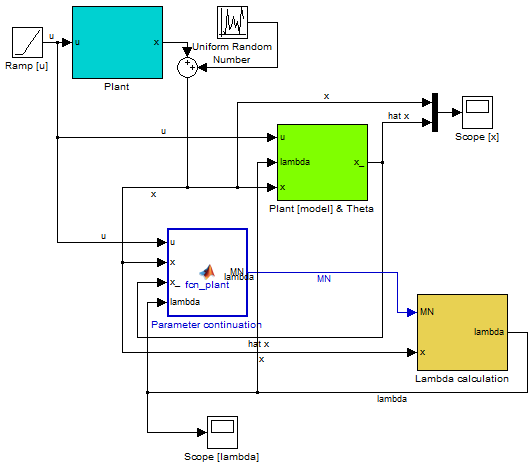}
\fi

Figure 9. Model of the system in Simulink.
\end{center}

\begin{center}
\ifpdf 
  \includegraphics[width=0.85\textwidth]{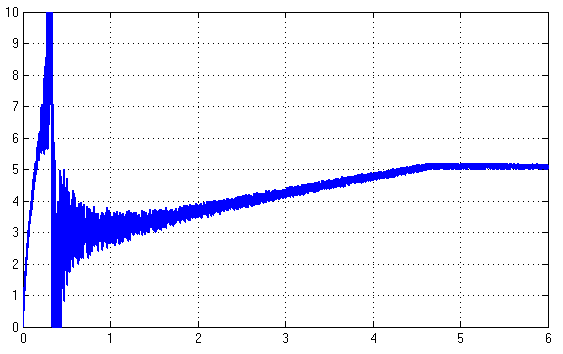}
\fi

Figure 10. Dynamics of the identified parameter $\hat \theta(t)$.
\end{center}

\begin{center}
\ifpdf 
  \includegraphics[width=0.85\textwidth]{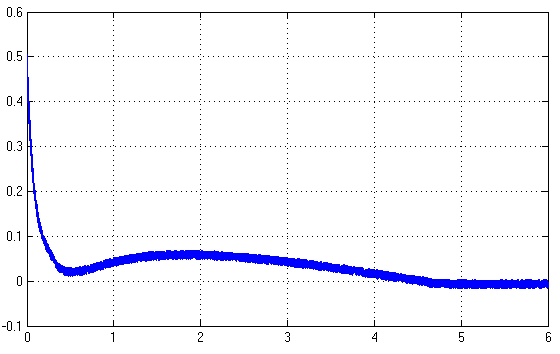}
\fi

Figure 11. Dynamics of the identification error $\hat x - x$.
\end{center}

\begin{center}
\ifpdf 
  \includegraphics[width=0.85\textwidth]{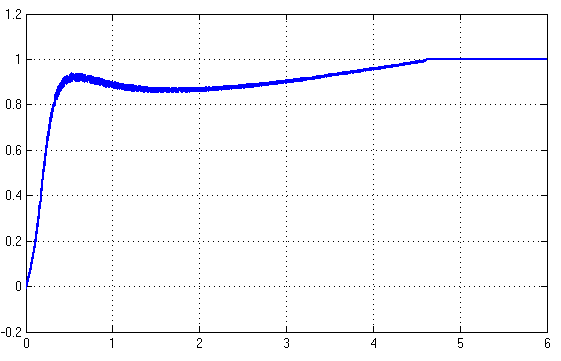}
\fi

Figure 12. Dynamics of parameter $\lambda$.
\end{center}

\subsection{Identification for the case of nonlinearly parametrized uncertainties}\label{sec:nonaffine_id}

In this section we consider the problem of system identification with nonlinear parametrization of uncertainty. A new interpretation of the parameter identification problem other than provided in section \ref{sec:affine_id} is given.

Let the following nonlinear system which uncertainties that are parametrized by a vector $\theta$ is given

\begin{equation}\label{plant_uncertain_}
\dot{x} = f(x,u,\theta)
\end{equation}

with initial conditions $x(0) = x_0$, known input signal $u(t)$ and with the measurable state $x$.

For the system \eqref{plant_uncertain_} the model is known, which uses a parameter estimate $\hat \theta$ instead of the true value $\theta$:

\begin{equation}\label{plant_uncertain_model}
\dot{\hat x} = f(\hat x, u, \hat \theta)
\end{equation}

The problem of the method described in \ref{sec:affine_id}, is that to zero identification error $\hat x - x$ the controllability of state $\hat x$ by parameter $\hat \theta$ is required in system \eqref{plant_uncertain_model}. In fact, adaptation algorithm can be restated as the tracking problem for \eqref{plant_uncertain_model} and input $\hat \theta$ with objective $\hat x \to x$.

Nevertheless, even if system with an output $\hat x$ is not completely controlled by $\hat \theta$ it is possible to solve the problem of parameter identification. The following simple idea is for this purpose offered: periodically in the discrete timepoints $t_i$ restart integration of model \eqref{plant_uncertain_model}, and for the initial state take the state of identified plant $\hat x(t_i) := x(t_i)$. This idea has its roots in the method of nonlinear model predictive control. Hence we conclude that during the time interval $t \in [t_i, t_{i+1})$ the state $\hat x(t)$ can be expressed by integral:

\begin{equation}\label{int_uncertain_model}
\begin{gathered}
\hat x(t) = \int_{t_i}^t f(\hat x(\tau),u(\tau),\hat \theta) d \tau \\
\hat x(t_i) = x(t_i)
\end{gathered}
\end{equation}

Thus, at the time $t_{i+1}$ before the next restart of the integration of \eqref{plant_uncertain_model} it is possible to define the following error signal

\begin{equation}\label{error_uncertain_model}
\begin{split}
e(\hat \theta,t_i) &= \int_{t_i}^{t_{i+1}} f(\hat x(\tau),u(\tau),\hat \theta) d \tau - x(t_{i+1}) \\
&= \int_{t_i}^{t_{i+1}} f(\hat x(\tau),u(\tau),\hat \theta) - f(x(\tau),u(\tau),\theta) d \tau
\end{split}
\end{equation}

Function $e(\hat \theta,t_i)$ defined in discrete moments of time $t_i$ and is subject to the minimization $e(\hat \theta,t_i) \to 0$ for $t_i \to \infty$.

Consider the continuous version of \eqref{error_uncertain_model} in order to understand the properties of the minimization problem for $e(\hat \theta,t_i)$. In the limit $t_i \to t_{i+1}$, it is possible to write an identification error as

\begin{equation}\label{error_uncertain_model_cont}
e(\hat \theta,t_i) \to \bar e(\hat \theta,t) = f(x(t),u(t),\hat \theta) - f(x(t),u(t),\theta)
\end{equation}

Should be noted a few features of the problem for solution of equation $\bar e(\hat \theta,t) = 0$:

- this problem is time-dependent (nonstationary), i.e. in general $\bar e(\hat \theta,t_1) \ne \bar e(\hat \theta,t_2)$ for $t_1 \ne t_2$,

- the position of the solution point $\theta$ is stationary, i.e. for $\hat \theta = \theta$ following $\bar e(\theta,t) = 0$ is true for all $t$.

The task of zeroing $e(\hat \theta,t) \to 0$ can be solved by a known method based on feedback linearization (Newton's method). Differentiating $e(\hat \theta,t)$ it is possible to obtain

\begin{equation}\label{derror_uncertain_model_cont}
\begin{split}
\dot{\bar e} &= \frac{\partial \bar e}{\partial \hat \theta} \cdot \dot {\hat \theta} + \frac{\partial \bar e}{\partial t} \\
&= \frac{\partial f}{\partial \hat \theta} \cdot \dot {\hat \theta} + \frac{\partial \bar e}{\partial t}
\end{split}
\end{equation}

Assuming that $\rank \dfrac{\partial f}{\partial \hat \theta} = n$, it is possible to define the following dynamics of the parameter estimation $\hat \theta$:

\begin{equation}\label{d_theta_uncertain}
\begin{split}
\dot{\hat \theta} &= \left ( \frac{\partial f}{\partial \hat \theta} \right )^{-1} \cdot v \\
&= - \left ( \frac{\partial f}{\partial \hat \theta} \right )^{-1} \cdot k \bar e
\end{split}
\end{equation}

where $k > 0$ is a feedback gain.

Substituting \eqref{d_theta_uncertain} in to \eqref{derror_uncertain_model_cont} we obtain the linearized dynamics

\begin{equation}\label{derror_uncertain_model_cl}
\begin{split}
\dot{\bar e} &= -k \cdot \bar e + \frac{\partial \bar e}{\partial t} \\
&= -k \cdot \bar e + \Delta
\end{split}
\end{equation}

where $\Delta = \dfrac{\partial \bar e}{\partial t}$.

\begin{Thm}\label{thm:conv} 

If 

\begin{equation}\label{rank_uncertain_cond}
\rank \dfrac{\partial f}{\partial \hat \theta} = n
\end{equation}

and 

\begin{equation}
| \Delta_i | < \infty
\end{equation}

it is always possible to choose a vector $k > 0$ in \eqref{d_theta_uncertain}, that neighborhood of a point $\hat \theta = \theta$ is globally asymptotically stable and convergence $\bar e \to \epsilon$ can be achieved for all initial $\hat \theta(0)$.

\end{Thm} 

\begin{Prf}
Outline of the proof can be described as follows. Because each component of the vector $\Delta = \partial \bar e/\partial t$ is limited it is always possible to choose constant $k > 0$ in such way that in some area $\mathcal{O}$ perturbation $\Delta$ does not affect the convergence of $\bar e$. Further analysis is carried out in complement $\bar{\mathcal {O}}$ of the set $\mathcal{O}$, which can be regarded as a neighborhood of $\hat \theta = \theta$. Expanding $\bar e$ in the neighborhood of $\hat \theta $ it can be shown that $\hat \theta = \theta$ is a locally stable.
\end{Prf}

\subsection{Numerical continuation for the identification of nonlinear parameterized uncertainties}\label{sec:nonaffine_id_alg}

Described in \eqref{sec:nonaffine_id} is the basis and the main idea of the method of identification. However, before formulation a practical adaptation algorithm with numerical parameter continuation we must develop the following aspects:

- discretization of adaptation law \eqref{d_theta_uncertain} for using the error signal $e(\hat \theta,t_i)$ from \eqref{error_uncertain_model},

- computation of a partial derivative $\partial e/\partial \hat \theta$ for discrete case,

- relaxing the rank condition \eqref{rank_uncertain_cond}, ie actual application of the parameter continuation methodology.

Denote the partial derivative $\partial e/\partial \hat \theta$ of function $e(\hat \theta,t_i)$, calculated at time $t_i$

\begin{equation}
D_{\hat \theta} e(\hat \theta,t_i) = \frac{\partial e(\hat \theta,t_i)}{\partial \hat \theta}
\end{equation}

From \eqref{error_uncertain_model} it is possible to write

\begin{equation}\label{de_dtheta_uncertain_model}
\begin{split}
D_{\hat \theta} e(\hat \theta,t_i) &= \frac{\partial}{\partial \hat \theta} \int_{t_i}^{t_{i+1}} f(\hat x(\tau),u(\tau),\hat \theta) d \tau - x(t_{i+1}) \\
&= \int_{t_i}^{t_{i+1}} \frac{\partial}{\partial \hat \theta} f(\hat x(\tau),u(\tau),\hat \theta) d \tau
\end{split}
\end{equation}

where the last integral is element-wise integration of the Jacobian matrix $f$.

The algorithm of the regulator is represented further continuously in time, discretization of signals is considered only for $\partial e/\partial \hat \theta$ and $e$, which are converted to continuous with zero-order hold transformation.

Hence \eqref{d_theta_uncertain} can be written at time interval $[t_i,t_{i+1})$

\begin{equation}
\begin{split}
\dot{\hat \theta} &= - k \left ( D_{\hat \theta} e(\hat \theta,t_i) \right )^{-1} \cdot e(\hat \theta,t_i) \\
&= -k \left ( \int_{t_i}^{t_{i+1}} \frac{\partial}{\partial \hat \theta} f(\hat x(\tau),u(\tau),\hat \theta) d \tau \right )^{-1} \cdot \\
& \quad \cdot \left ( \int_{t_i}^{t_{i+1}} f(\hat x(\tau),u(\tau),\hat \theta) d \tau - x(t_{i+1}) \right )
\end{split}
\end{equation}

Let's consider the remained subtask, namely application of parameter continuation strategy to the problem of adaptive identification in a case when the rank condition \eqref{rank_uncertain_cond} isn't satisfied. As mentioned above, the identification problem is reduced to a nonlinear equation $e(\hat \theta,t_i) = 0$ with disturbances. Therefore, the material of the theorem \ref{Thm:sol} can be directly applied to the system $y = e(\hat \theta,t_i)$, $\dot {\hat \theta} = u$, $\mathcal{B} = 0$ and $r_i = 1$, and considering regulation of the convergence rate using \eqref{Sol_gamma}.

Let's construct the homotopy mapping

\begin{equation}\label{H_uncertain}
H = \lambda \cdot e + (1 - \lambda) \cdot (\hat \theta - \hat \theta_0)
\end{equation}

Further it is possible write the matrix $\mathcal{A}$
\begin{equation}\label{A_uncertain}
\mathcal{A} = \left( \lambda D_{\hat \theta} e(\hat \theta,t_i) + (1-\lambda) E ~~~ e(\hat \theta,t_i) - \hat \theta + \hat \theta_0 \right )
\end{equation}

where $E$ is unit diagonal matrix.

Finally, we obtain the dynamics of the solutions as a special case of \eqref{Sol_gamma}

\begin{equation}\label{Sol_uncertain_gamma}
\begin{pmatrix} \dot{\hat \theta} \\ \dot \lambda \end{pmatrix} = \alpha \cdot \gamma \cdot \tau(\mathcal{A}) - k Q (\mathcal{A} \cdot Q)^+ H,
\end{equation}

where coefficient $\alpha > 0$ determines the priority of change $\lambda$, coefficient $\gamma > 0$ controls the rate of convergence and $k > 0$ is a proportional feedback gain, $Q = diag(1,...1,\gamma)$ is square diagonal matrix, the last diagonal element of which is equal to $\gamma$.

\subsection{Illustrative example for the system identification with nonlinear parametrization}

For example, consider the following nonlinear oscillation system

\begin{equation}\label{plant_ex_uncertain}
\begin{gathered}
\dot x_1 = x_2 \\
\dot x_2 = -\omega^2 \sin(\theta  x_1) + u
\end{gathered}
\end{equation}

performs forced oscillations under the action of disturbance $u(t) = \sin(4 t)$. Initial conditions are given $x(0) = (0,0)^T$.

The identified parameter is $\theta = 0.75$. Frequency value $\omega$ is known for controller: $\omega = 2$.

Consider two solutions: 

- abstract solution under the assumption that information of the derivatives $f(x(t),u(t),\theta)$ is available,

- complete solution using only observable states and signals $D_{\hat \theta} e(\hat \theta,t_i)$, $e(\hat \theta,t_i)$.

The first solution, discussed below, is only illustrative and can not be implemented because the signal $\dot x$ (ie. $f(x(t),u(t),\theta)$) in practice is not available for the measurement. However, this problem can be solved directly without time discretization.

Let's associate with \eqref{plant_ex_uncertain} model of the form

\begin{equation}\label{model_ex_uncertain}
\begin{gathered}
\dot {\hat x}_1 = \hat x_2 \\
\dot {\hat x}_2 = -\omega^2 \sin(\hat \theta {\hat x_1}) + u
\end{gathered}
\end{equation}

with the same as for \eqref{plant_ex_uncertain} initial conditions $\hat x(0) = (0,0)^T$.

Parameter identification implemented as parameter continuation for zeroing the error signal $\bar e(\hat \theta,t)$, which is considered as stationary function from $\hat \theta$ with perturbations.

Block diagram of the model in MATLAB/Simulink is shown in Figure 13. The model consists of the following subsystems: 

- "Plant" is an identified system with unknown parameter $\theta$ modeled by \eqref{plant_ex_uncertain}, 

- "Plant Model" is model of identified object for computation of $\dot {\hat x}$ via \eqref{model_ex_uncertain}, 

- "Derivative calc" is subsystem for computation of $\partial \dot \hat {x_2} / \partial \hat \theta$, 

- "Feedback linearization \& continuation" is embedded MATLAB function, that implements the parameter continuation dynamics \eqref{Sol_uncertain_gamma}.

Note that in the continuous case, while minimizing the difference of derivatives $\bar e(\hat \theta,t)$ in \eqref{model_ex_uncertain} used $\hat x := x$.

For the parameters of the matrix $\mathcal{A}$ used in \eqref{Sol_uncertain_gamma} is selected: $\gamma = 0.1$, $k = 10$, $\alpha = 1$.

The simulation results are presented in Figures 14-16.

\begin{center}
\ifpdf 
  \includegraphics[width=1.0\textwidth]{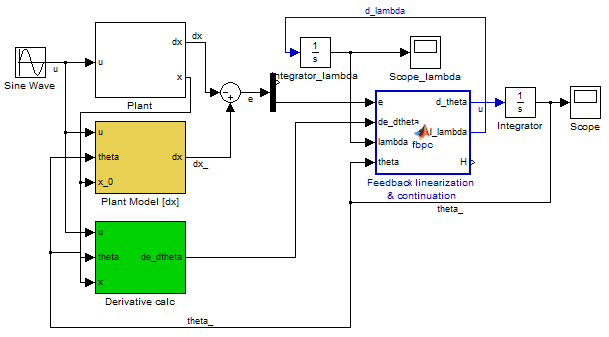}
\fi

Figure 13. The model in Simulink.
\end{center}

\begin{center}
\ifpdf 
  \includegraphics[width=0.85\textwidth]{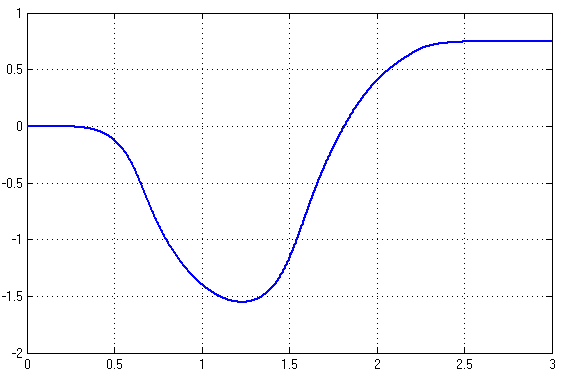}
\fi

Figure 14. Dynamics of the identified parameter $\hat \theta(t)$.
\end{center}

\begin{center}
\ifpdf 
  \includegraphics[width=0.85\textwidth]{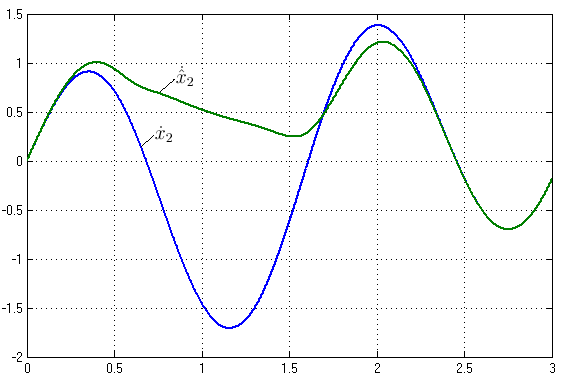}
\fi

Figure 15. Derivative $\dot x_2$ of the plant state variable $x_2$ and derivative $\dot {\hat x}_2$ for model state variable.
\end{center}

\begin{center}
\ifpdf 
  \includegraphics[width=0.85\textwidth]{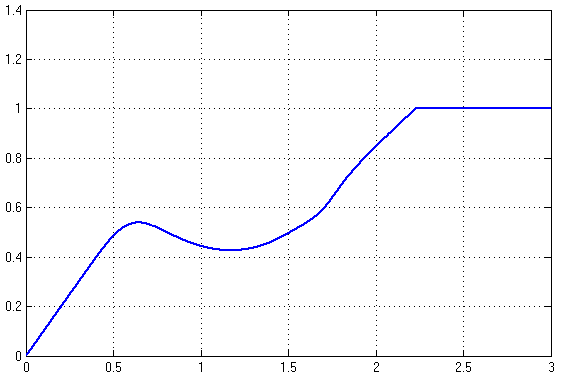}
\fi

Figure 16. Dynamics of parameter $\lambda(t)$.
\end{center}

From the obtained data it is evident that the algorithm converges to the desired exact value $\hat \theta = 0.75$ of estimated parameter.

The second solution, which uses only the signals $D_{\hat \theta} e(\hat \theta,t_i)$ and $e(\hat \theta,t_i)$ based on the method described in section \ref{sec:nonaffine_id_alg} and is a primary interest for the implementation.

Block diagram of the model in MATLAB/Simulink is shown in Figure 13. The model consists of the following subsystems: 

- "Plant" is an identified system with unknown parameter $\theta$ modeled by \eqref{plant_ex_uncertain}, 

- "Plant Model" is model of identified object for computation of  $e(\hat \theta,t_i)$ via \eqref{error_uncertain_model}, 

- "Derivative calc" is subsystem for computation of $D_{\hat \theta} e(\hat \theta,t_i)$ by \eqref{de_dtheta_uncertain_model}, 

- "Feedback linearization \& continuation" is embedded MATLAB function, that implements the parameter continuation dynamics \eqref{Sol_uncertain_gamma}.

Sampling time $\Delta t = t_{i+1} - t_{i}$ is selected 0.1 sec. 

The subsystem "Plant Model" implement integration of \eqref{model_ex_uncertain} and periodic updating of the initial state $\hat x(t_i) = x(t_i)$ after an interval $\Delta t$, thus realizing the equation \eqref{error_uncertain_model}. 

Subsystem "Derivative calc" calculates the integral of the partial derivative $\partial \dot {\hat x}_2 / \partial \theta {\hat x_1}$

\begin{equation}
D_{\hat \theta} e(\hat \theta,t_i) = -\omega^2 \int_{t_i}^{t_{i+1}} \cos(\hat \theta \cdot x_1) \cdot x_1 d\tau
\end{equation}

For the parameters of the matrix $\mathcal{A}$ used in \eqref{Sol_uncertain_gamma} is selected: $\gamma = 0.05$, $k = 10$, $\alpha = 0.5$.

The simulation results are presented in Figures 18-20.

\begin{center}
\ifpdf 
  \includegraphics[width=1.0\textwidth]{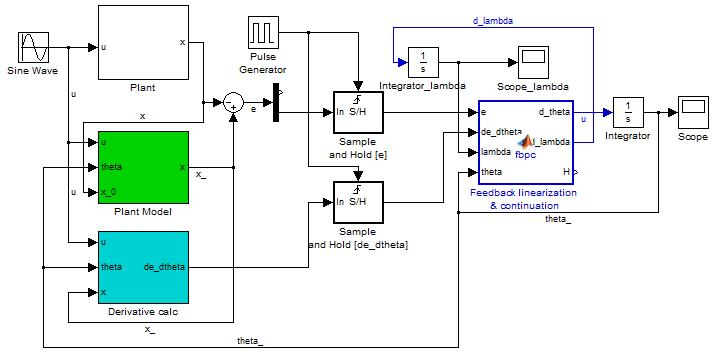}
\fi

Figure 17. The model in Simulink.
\end{center}

\begin{center}
\ifpdf 
  \includegraphics[width=0.85\textwidth]{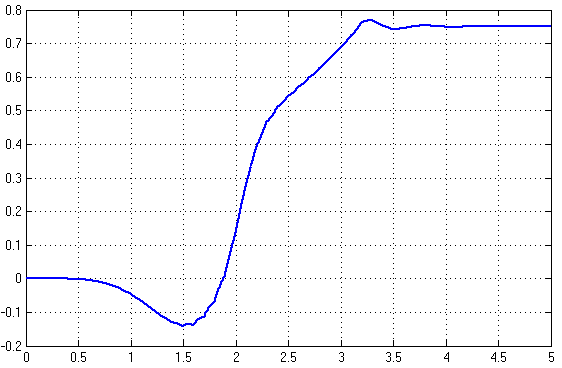}
\fi

Figure 18. Dynamics of the identified parameter $\hat \theta(t)$.
\end{center}

\begin{center}
\ifpdf 
  \includegraphics[width=0.85\textwidth]{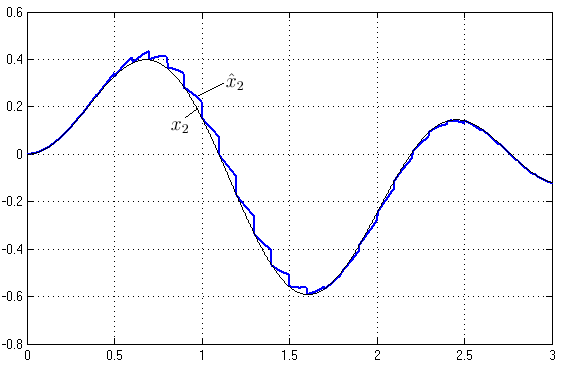}
\fi

Figure 19. The state $x_2$ of identified plant and model $\hat x_2$.
\end{center}

\begin{center}
\ifpdf 
  \includegraphics[width=0.85\textwidth]{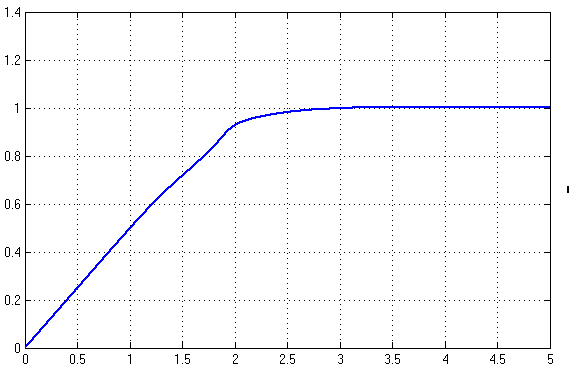}
\fi

Figure 20. Dynamics of parameter $\lambda(t)$.
\end{center}

From the obtained data it is evident that the algorithm converges to the desired value of the estimated parameter $\hat \theta = 0.75$. The dynamics of $\lambda(t)$ not experiencing the bifurcation bends which is caused by discretization of model.

\subsection{Adaptive control}

In this section we only point out that on-line parameter identification is central to a particular type of adaptive control -- Model Identification Adaptive Controller (MIAC) \cite{20}.

The aim of adaptive control with reference model (model reference adaptive control, MRAC) is to achieve closed-loop behavior prescribed by dedicated model. In the case of MIAC type control the continuous adjustment of the model parameters is performed to minimize error between the model and the plant. So there is no need for a reference model describing the desired characteristics of a closed-loop.

\begin{center}
\ifpdf 
  \includegraphics[width=0.85\textwidth]{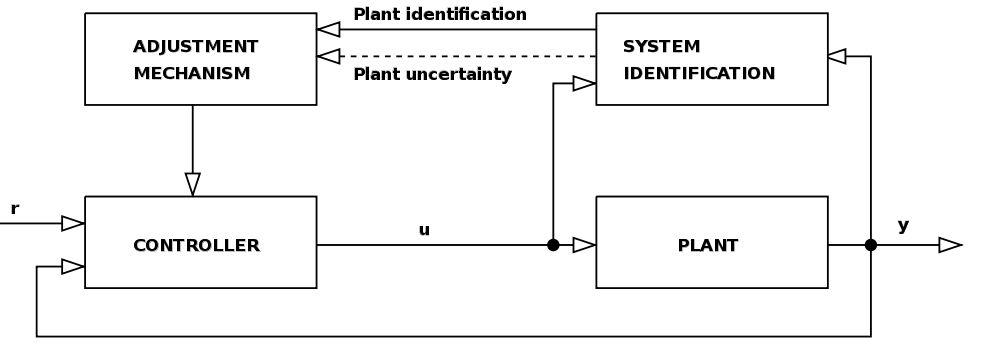}
\fi

Figure 21. Adaptive control structure with the identification of parameters (from \cite{20}).
\end{center}

Thus, both at a output regulation stage and at a stage of parameter identification the application of parameter continuation methodology as a part of the adaptive regulator is possible.

\section*{Conclusion}

A survey is given of the parameter continuation method used for the control and parameter identification of nonlinear systems. Development of the computational basis for the implementation of nonlinear systems allow to consider more sophisticated regulator algorithms that can handle more difficult to control systems (in particular, failed to feedback linearization).

Without claiming to cover all of the possible applications of parameter continuation methods in the paper following control problems have been considered:

- output setpoint control for the affine nonlinear systems,

- output setpoint control for the nonaffine nonlinear systems in general form,

- online parameter identification of nonlinear affine systems,

- online parameter identification of nonlinear nonaffine plants.

Further work will be carried out in two directions:

1. First, the practical meaning of the proposed methods has yet to be confirmed by the application of this algorithm to real existing systems, the control of which was difficult before. In our view the application to systems with multiple equilibrium states between which the phase transition occurs \cite{21} is quite promising.

2. As part of the general methodology of parameter continuation, namely parameterized transition from simple to complex problems, the offered methods shall receive more abstract theoretical basis. Promising fields of application is category theory and categorical homotopy theory to the formulation of the proposed algorithms and related phenomenons in the control theory \cite{22}.

\ifpdf 
In conclusion, it should be added that all Simulink models considered as illustrations for this paper can be downloaded here \url{https://sites.google.com/site/akpc806a/ParameterContinuationFiles.rar}
\fi

\end{document}